\ifx\shlhetal\undefinedcontrolsequence\let\shlhetal\relax\fi

\documentstyle[12pt]{article}

\newtheorem{theorem}{\sc Theorem}[section]

\newtheorem{proposition}[theorem]{\sc Proposition}
\newtheorem{lemma}[theorem]{\sc Lemma}
\newtheorem{corollary}[theorem]{\sc Corollary}
\newtheorem{claim}{\sc Claim}[theorem]

\newcommand{\proof}{\noindent {\sc Proof. }} 

\newcounter{defi}
\newenvironment{defi}{\stepcounter{defi}\medskip\noindent {\sc
Definition \thesection.\thedefi\ }}{\par\medskip} 

\newcommand{\open}{\Bbb}

\newcommand{\se}{\subseteq}
\newcommand{\set}[2]{\{#1 \colon #2\}} 
\newcommand{\fin}{$\Box$\medskip} 
\def\deq{\mathop=\limits^{\rm def}} 

\newcommand{\ga}{\alpha}
\newcommand{\gb}{\beta}
\newcommand{\gd}{\delta}
\newcommand{\gk}{\kappa}
\newcommand{\gl}{\lambda}
\newcommand{\gm}{\mu}
\newcommand{\gf}{\varphi}
\newcommand{\gth}{\theta}
\newcommand{\go}{\omega}

\newcommand{\gG}{\Gamma}
\newcommand{\gD}{\Delta}

\newcommand{\gL}{\Lambda}

\newcommand{\ha}{\aleph}

\def\ltc{\mathord{<}}

\def\rest{\mathord{\restriction}}
\newcommand{\force}{\Vdash}

\newcommand{\He}[1]{\mbox{\rm H}(#1), \in}

\newcommand{\enf}{essentially non-free}
\renewcommand{\iff}{if and only if}

\newcommand{\isub}[1]{\prec_{\infty #1}}
\newcommand{\cp}[1]{$\mbox{\rm CP}_{#1}$}
\newcommand{\lsub}[1]{$\mbox{\rm L}_{\infty #1}$}

\input mssymb
\title{Almost Free Algebras}
\author{Alan H. Mekler\thanks{Research supported by NSERC grant
\#9848}\\
 Department of Mathematics and Statistics \\ 
 Simon Fraser University \\
 Burnaby, B.C. V5A 1S6 
CANADA  \and Saharon Shelah\thanks{Publication
\#366. Research supported by the BSF}\\ Institute of Mathematics\\ The
Hebrew University\\Jerusalem, ISRAEL\\	and\\ Department of Mathematics\\
Rutgers University\\New Brunswick, New Jersey, USA}
\begin{document}
\maketitle
\section{Introduction}

In this paper we adopt the terminology of universal algebra.  So by a
free algebra we will mean that a variety (i.e., an equationally defined
class of algebras) is given and the algebra is free in that variety. We
will also assume that the language of any variety is countable.

In this paper the investigation of the almost free algebras is continued.
An algebra is said to be almost free if ``most'' of its subalgebras of
smaller cardinality are free. For some varieties, such as groups and
abelian groups, every subalgebra of a free algebra is free. In those
cases ``most'' is synonomous with ``all''.  In general there are several
choices for the definition of ``most''. In the singular case, if the
notion of ``most'' is strong enough, then any almost free algebra of
singular cardinality is free \cite{Sing}. So we can concentrate on the
regular case. In the regular case we will adopt the following definition.
If $\gk$ is a regular uncountable cardinal and $A$ is an algebra of
cardinality $\gk$ then $A$ is {\em almost free} if there is a sequence
$(A_\ga\colon \ga < \gk)$ of free subalgebras of $A$ such that: for all
$\ga$, $|A_\ga| < \gk$; if $\ga < \gb$, then $A_\ga \se A_\gb$; and if
$\gd$ is a limit ordinal then $A_\gd = \bigcup_{\ga < \gd} A_\ga$. (In
\cite{EMB}, such a chain is called a $\gk${\em -filtration}.) It should
be noted that this definition is not the same as the definition in
\cite{Sing}. The definition there is sensitive to the truth of Chang's
conjecture (see the discussion in
\cite{EMB} Notes to Chapter IV).

There are two sorts of almost free algebras; those which are essentially
free and those which are essentially non-free. An algebra $A$ is {\em
essentially free} if $A \ast F$ is free for some free algebra $F$.  Here
$\ast$ denotes the free product. (Since our algebras will always be
countably free the free product is well defined.) For example, in the
variety of abelian groups of exponents six, a free group is a direct sum
of copies of cyclic groups of order 6. The group $\bigoplus_{\ha_1}C_3
\oplus \bigoplus_{\ha_0}C_2$ is an almost free algebra of cardinality
$\ha_1$ which is not free but is essentially free. An algebra which is
not essentially free is {\em essentially non-free}.  In \cite{EM}, the
construction principle, abbreviated CP, is defined and it is shown that
for any variety there is an essentially non-free almost free algebra of
some cardinality if and only if there is an \enf\ almost free algebra of
cardinality $\aleph_1$ \iff\ the construction principle holds in that
variety. As well if V = L holds then each of the above equivalents is
also equivalent to the existence of an \enf\ almost free algebra in all
non-weakly compact regular cardinalities.

In this paper we will investigate the \enf\ spectrum of a variety.  The {\em
\enf\ spectrum} is the class of uncountable cardinals $\gk$ in which there
is an \enf\ algebra of cardinality $\gk$ which is almost free. This class
consists entirely of regular cardinals (\cite{Sing}). In L, the \enf\
spectrum of a variety is entirely determined by whether or not the
construction principle holds. As we shall see the situation in ZFC may be
more complicated.

For some varieties, such as groups, abelian groups or any variety of
modules over a non-left perfect ring, the \enf\ spectrum contains not
only $\aleph_1$ but $\aleph_n$ for all $n > 0$.  The reason for this
being true in ZFC (rather than under some special set theoretic
hypotheses) is that these varieties satisfy stronger versions of the
construction principle. We conjecture that the hierarchy of construction
principles is strict, i.e., that for each $n > 0$ there is a variety
which satisfies the $n$-construction principle but not the $n +
1$-construction principle. In this paper we will show that the
$1$-construction principle does not imply the $2$-construction principle.

After these examples are given there still remains the question of
whether these principles actually reflect the reason that there are \enf\
$\gk$-free algebras of cardinality $\gk$. Of course, we can not hope to
prove a theorem in ZFC, because of the situation in L (or more generally
if there is a non-reflecting stationary subset of every regular
non-weakly compact cardinal which consists of ordinals of cofinality
$\go$).  However we will prove that, assuming the consistency of some
large cardinal hypothesis, it is consistent that a variety has an
\enf\ almost free algebra of cardinality $\aleph_n$
\iff\ it satisfies the $n$-construction principle. (We will also show
under milder hypotheses that it is consistent that the various classes
are separated.) 

\begin{defi} 
A variety ${\cal V}$ of algebras satisfies the $n$-construction
principle, \cp{n}, if there are countably generated free algebras $H
\subseteq I \subseteq L$ and a partition of $\omega $ into $n$ infinite
blocks (i.e. sets) $s^1, \ldots, s^n$ so that
\begin{quotation}
\noindent
(1) $H$ is freely generated by $\{h_m\colon m < \omega \}$, and for every
subset $J \subseteq \omega $ if for some $k$, $J \cap s^k$ is finite then
the algebra generated by $\{h_m\colon m \in J\}$ is a free factor of $L$;
and
\newline (2) $L = I \ast F(\omega )$ and $H$ is not a free factor of
$L$.
\end{quotation} 
\end{defi}
 
Here $F(\omega)$ is the free algebra on $\aleph _0$ generators, and $H$
is a free factor of $L$, denoted $H | L$, means that there is a free
algebra $G$ so that $H \ast G = L$.

The construction principle of \cite{EM} is the principle we have called
\cp{1}. The known constructions of an almost free algebra from \cp{n}
seem to require the following set theoretic principle. (The definition
that follows may be easier to understand if the reader keeps in mind that
a $\gl$-system is a generalization of a stationary set consisting of
ordinals of cofinality $\go$.)

\begin{defi} (1) A $\lambda$-{\em set} of height $n$ is a
subtree $S$ of $^{\leq n }\!\lambda $ together with a cardinal $\lambda
_\eta $ for every $\eta \in S$ such that $\lambda _\emptyset = \lambda $,
and:
 
(a) for all $\eta \in S$, $\eta $ is a final node of $S$ if and only if
$\lambda _\eta = \aleph _0;$

(b) if $\eta \in S \setminus S_f{ }$, then $\eta \hat{\ }\langle
\beta
\rangle \in S$ implies $\beta \in \lambda _\eta $, $\lambda _{\eta
\hat{\ }\langle \beta \rangle } < \lambda _\eta $ and $E_\eta \deq \{\beta
< \lambda _\eta \colon \eta \hat{\ }\langle \beta \rangle \in S\}$ is
stationary in $\lambda _\eta $ (where $S_f\stackrel{{\rm def}}{=}\{\eta\in
S: \eta\mbox{ is }\triangleleft\mbox{-maximal in } S\}$ is the family of
final nodes of $S$).
 
\medskip
 
(2) A $\lambda ${\em -system} of height $n$ is a $\lambda$-set of height
$n$ together with a set $B_\eta $ for each $\eta \in S$ such that
$B_\emptyset = \emptyset $, and for all $\eta \in S \setminus S_f$:
 
(a) for all $\beta \in E_\eta $, $\lambda _{\eta \hat{\ }\langle
\beta \rangle } \leq |B_{\eta \hat{\ }\langle \beta \rangle }| <
\lambda _\eta ;$
 
(b) $\{B_{\eta \hat{\ }\langle \beta \rangle }\colon \beta \in E_\eta
\}$ is an increasing continuous chain of sets, i.e., if $\beta < \beta'$
are in $E_\eta $, then $B_{\eta \hat{\ }\langle \beta \rangle}
\subseteq B_{\eta\hat{\ }\langle \beta '\rangle }$; and if $\sigma $ is a
limit point of $E_\eta $ (i.e. $\sigma=\sup(\sigma\cap E_\eta)\in E_\eta$),
then $B_{\eta \hat{\ }\langle \sigma \rangle } = \cup\{B_{\eta\hat{\
}\langle \beta \rangle }\colon \beta < \sigma $, $\beta \in E_\eta \}$.   
 
\medskip
 
(4) For any $\lambda $-system $\Lambda = (S, \lambda _\eta , B_\eta
\colon \eta \in S)$, and any $\eta \in S$, let $\bar B_\eta = \cup
\{B_{\eta \rest m}\colon m \leq \ell (\eta )\}$.  Say that a family
${\cal S}$ of countable sets is {\em based on $\Lambda $} if ${\cal S}$ is
indexed by $S_f$, and for every $\eta \in S_f$, $s_\eta
\subseteq \bar B_\eta $. 
 
\medskip

A family ${\cal S}$ of countable sets is {\em free} if there is a
transversal of $\cal S$, i.e., a one-one function $f$ from $\cal S$ to $\cup
\cal S$ so that for all $s \in \cal S$ $f(s) \in s$. A family of countable
sets is, {\em almost free} if every subfamily of lesser cardinality is free.
\end{defi}

Shelah, \cite{Reg}, showed that the existence of an almost free abelian
group of cardinality $\gk$ is equivalent to the existence of an almost free
family of countable sets of cardinality $\gk$. The proof goes through
$\gl$-systems.  In \cite{EMB}, the following theorem is proved (although not
explicitly stated, see the proof of theorem VII.3A.13).

\begin{theorem}
\label{existence}
If a variety satisfies \cp{n} and $\gl$ is a regular cardinal such that
there is a $\gl$-system, $\gL$, of height $n$ and an almost free family of
countable sets based on $\gL$, then there is an \enf\ algebra of cardinality
$\gl$ which is almost free.
\end{theorem}

\medskip
\noindent{\sc Conjecture:}\ \ \ The converse of the theorem above is true.
I.e. for each regular cardinal $\lambda>\aleph_0$ and every variety the
following two conditions are equivalent:
\begin{description}
\item[$(\alpha)$\ \ \ ] for some $n<\omega$ the variety satisfies the
principle \cp{n} and there exists a $\lambda$-system $\gL$ of height $n$ and
an almost free family of countable sets based on $\gL$
\item[$(\beta)$\ \ \ ] there exists an essentially non free almost free
algebra (for the variety) of cardinality $\lambda$. 
\end{description}
Since $\{n: \mbox{ the variety satisfies \cp{n}}\}$ is an initial segment of
$\omega$ we could conclude that it is a theorem of ZFC that there are at
most $\ha_0$ \enf\ spectra.

\medskip

Although we will not discuss essentially free algebras in this paper, these
algebras can be profitably investigated. The {\em essentially free spectrum}
of a variety is defined as the set of cardinals $\gk$ so that there is an
almost free non-free algebra of cardinality $\gk$ which is essentially free.
The conjecture is that the essentially free spectrum of a variety is either
empty or consists of the class of successor cardinals. For those varieties
for which \cp{1} does not hold, i.e., the \enf\ spectrum is empty, the
conjecture is true \cite{linf}. It is always true that the essentially free
spectrum of a variety is contained in the class of successor cardinals. (A
paper which essentially verifies it is in preparation.)

A notion related to being almost free is being $\gk$-free where $\gk$ is an
uncountable cardinal. An algebra is $\gk$-free if ``most'' subalgebras of
cardinality less than $\gk$ are free. There are various choices for the
definition of ``most'' and the relations among them are not clear. For a
regular cardinal $\gk$ we will say that $A$ is {\em$\gk$-free} if there is a
closed unbounded set in ${\cal P}_\gk(A)$ (the set of subsets of $A$ of
cardinality less than $\gk$) consisting of free algebras. Note that an
almost free algebra of cardinality $\gk$ is $\gk$-free.  One important
associated notion is that of being \lsub{\gk}-free; i.e., being
\lsub{\gk}-equivalent to a free algebra. A basic theorem is that:

\begin{theorem}
If an algebra is $\gk^+$-free, then it is \lsub{\gk}-free.
\end{theorem}

\proof See  \cite{Sing} 2.6(B). (Note $\gk^+$-free 
in the sense here implies $E^{\gk^+}_{\gk^+}$-free as defined there.)
\fin

We will use some of the notions associated with the
\lsub{\gk}-free algebras. Suppose $\gk$ is a  cardinal and $A$ is an
algebra (in some fixed variety). A subalgebra $B$ which is
$\ltc\gk$-generated is said to be {\em $\gk$-pure} if Player I has a
winning strategy in the following game of length $\go$.
\begin{quote}
Players I and II alternately choose an increasing chain $B = B_0 \se B_1
\se \ldots \se B_n \se \ldots $ of subalgebras of A each of which is
$\ltc\gk$-generated. Player I wins a play of the game if for all $n$,
$B_{2n}$ is free and $B_{2n}$ is a free factor of $B_{2n +2}$.
\end{quote}
If $B\subseteq A$ is not $<\kappa$-generated (but $B\subseteq A$ are free)
we just ask $B_{n+1}$ to be $<\kappa$-generated over $B_n$ for each $n$
(used e.g. in the proof of \ref{add1-thm}). 

The choice of the term $\gk$-pure is taken from abelian group theory.  The
following theorem sums up various useful facts. Some of the results are
obvious others are taken from \cite{Kuek}.

\begin{theorem}\ 
\begin{trivlist}
\item[\rm (1)] An algebra is \lsub{\gk}-free \iff\ every subset of
cardinality less than $\gk$ is contained in a $\ltc\gk$-generated algebra
which is $\gk$-pure.
\item[\rm (2)] If $F$ is a free algebra then a subalgebra is
$\gk$-pure \iff\ it is a free factor.
\item[\rm (3)] In any \lsub{\gk}-free algebra the set of $\gk$-pure
subalgebras is $\gk$-directed under the relation of being a free factor.
\item[\rm (4)] If $\gk < \gl$ and $A$ is \lsub{\gl}-free, then any
$\gk$-pure subalgebra is also $\gl$-pure.
\end{trivlist}
\end{theorem}

Notice that part (2) of the theorem above implies that for $\gk$, if $A$ is
\lsub{\gk}-free then there is a formula of \lsub{\gk} which defines the
$\gk$-pure subalgebras (of $A$, but the formula depends on $A$).

We will use elementary submodels of appropriate set theoretic universes on
many occasions. We say that a cardinal $\chi$ {\em is large enough} if
$(\He{\chi})$ contains as elements everything which we are discussing. If
$A$ and $B$ are free algebras which are subalgebras of some third algebra
$C$, then by $A+ B$, we denote the algebra generated by $A \cup B$ and
define $B/A$ to be {\em free} if any (equivalently, some) free basis of $A$
can be extended to a free basis of $A + B$.  Similarly for $\gk$ a regular
uncountable cardinal, if $A + B$ is $\gk$-generated over $A$ then we say
that $B/A$ is {\em almost free} if there is a sequence $(B_\ga\colon\ga <
\gk)$ so that: for all $\ga$, $|B_\ga| < \gk$; for all $\ga < \gb$, $B_\ga
\se B_\gb$; if $\gd$ is a limit ordinal then $B_{\gd} = \bigcup_{\ga < \gd}
B_\ga$; $A + B = A +\bigcup_{\ga < \gk} B_\ga$; and for all $\ga$, $B_\ga/A$
is free. The notions of {\em \enf, strongly $\gk$-free} etc.\ for pairs are
defined analogously. The following lemma is useful.

\begin{lemma}
\label{sub-lemma}
Suppose $A \se B$ are free algebras and $N \prec (\He{\chi})$, where $\chi$
is large enough. If $A, B \in N$ and $(B \cap N)/(A \cap N)$ is free then
$(B \cap N)/A$ is free.
\end{lemma}

\proof Let $X \in N$ be a free basis of $A$. By elementariness, $A
\cap N$ is freely generated by $X \cap N$. Choose $Y$ so that $(X \cap
N) \cup Y$ is a free basis for $B \cap N$. We now claim that $X \cup Y$
freely generate the algebra they generate (namely $A + (B \cap N)$).
Suppose not. Then there are finite sets $Y_1 \se Y$ and $X_1 \se X$ such
that $Y_1 \cup X_1$ satisfy an equation which is not a law of the variety.
By elementariness, we can find $X_2 \se X\cap N$ so that $Y_1
\cup X_2$ satisfies the same equation. This is a contradiction. \fin

\section{\cp{1} does not imply \cp{2}}
\setcounter{theorem}{0}
\setcounter{defi}{0}

In this section we will present an example of a variety which satisfies
\cp{1} but not \cp{2}. The strategy for producing the example is quite
simple. We write down laws which say that the variety we are defining
satisfies \cp{1} and then prove that it does not satisfy \cp{2}. We believe
the same strategy will work for getting an example which satisfies \cp{n}
but not \cp{n+1}. However there are features in the proof that the strategy
works for the case $n = 1$ which do not generalize.

The variety we will build will be generated by projection algebras.  A {\em
projection algebra} is an algebra in which all the functions are projections
on some coordinate.  If a variety is generated by (a set of) projection
algebras, then it is not necessarily true that every algebra in the variety
is a projection algebra. For example, there may be a binary fuction $f$
which in one algebra is projection on the first coordinate and in another is
projection on the second coordinate.

In a variety generated by projection algebras there is a very simple
characterization of the free algebras. It is standard that a free algebra in
a variety is a subalgebra of a direct product of generators of the algebra
which is generated by tuples so that for any equation between terms there
is, if possible, a coordinate in which the equation fails for the tuples
(see for example Theorem~11.11 of \cite{BS}). In a variety generated by
projection algebras the free algebra on $\gk$ generators is the subalgebra
of the product of the various projection algebras on $\gk$ generators which
is generated by $\gk$ elements which differ pairwise in each coordinate.

\begin{theorem} 
There is a variety satisfying \cp{1} but not \cp{2}.
\end{theorem}

\proof  To begin we fix various sets of constant symbols: $\{c_{m,n}
\colon m, n < \go \}$ and $\{d_n \colon n < \go\}$. The intention is to
define an algebra $I$ such that for all $m$, $\{c_{m,n} \colon n < \go \}\cup
\{d_n \colon n < m \}$ will be a set of free generators for an algebra in
our variety. In particular, $\{c_{0n} \colon n < \go \}$ will be a set of
free generators. We intend that $H$ will be the algebra generated by $\{d_n
\colon n < \go \}$ and $I$ will be the whole algebra (in the definition of
\cp{1}). We define the language and some equations by induction. We have to
add enough function symbols so that for each $m$,
$\{c_{m,n}:n,\omega\}\cup\{d_n:n<m\}$ generates the whole $I$, it suffices
that it generates all $c\in\{c_{k,n}:k,n<\omega\}\cup\{d_n:n<\omega\}$.
However, while doing this we have still to make each
$\{c_{m,n}:n<\omega\}\cup\{d_n:n<m\}$ free. At each stage we will add a new
function symbol to the language and consider a pair consisting of a constant
symbol and a natural number $m$. (Note that the constants are not in the
language of the variety.) The enumeration of the pairs should be done in
such a way that each pair consisting a constant symbol and a natural number
is enumerated at some step. Since this is a routine enumeration we will not
comment on it, but assume our enumeration has this property.  Also at each
stage we will commit ourselves to an equation.
   
For the remainder of this proof we will let the index of the constant
$d_n$ be $n$ and the index of the constant $c_{m,n}$ be $m+n$. At stage
$n$ we are given a constant $t_n$ (so that $t_n \in
\{c_{m,k}, d_k\colon k, m \in \go\}$) and a natural number $m_n$, we
now add a new function symbol $f_n$ to the language where the arity of $f_n$
is chosen to be greater than $m_n$ plus the sum of the indices of $t_n$ and
all the constant symbols which which have appeared in the previous
equations. (No great care has to be taken in the choice of the arity, it
just has to increase quickly.)  Now we commit ourselves to the new equation
$$t_n = f_n(d_i (i < m_n ), c_{m_n,j} (j < k_n))\leqno(*)$$
where the arity of $f_n$ is $m_n + k_n$.

The variety we want to construct has vocabulary $\tau$, the set of function
symbols we introduced above. We use a subsidiary vocabulary $\tau'$ which is
$\tau \cup \{c_{m, n}, d_n\colon m, n <\go\}$. Let $K_0$ be the family of
$\tau'$-algebras which are projection algebras satisfying the equations
$(*)$ whose universe consists of $\{ c_{0,n} \colon n <
\go \}$ such that for all $m$, the interpretations of $c_{m,n}$ $( n <
\go)$ and $d_n$ $( n < m)$ are pairwise distinct. Let ${ \bf K }$ be the
class of $\tau$-reducts of members of $K_0$.  We will shortly prove that
${ \bf K }$ is non-empty. If we assume this for the moment, then it is
clear that the variety generated by ${ \bf K }$ satisfies (1) in the
definition of \cp{1} with $ \{d_n \colon n < \omega \}$ standing for $
\{h_n \colon n < \omega\}$. More exactly in the direct product of the
elements of ${ \bf K }$, for all $m$, $\{c_{m,n}\colon n < \go \} \cup
\{d_n \colon n < m \}$ freely generates a subalgebra. The choice of
equations guarantee that all the subalgebras are the same.  The proof
that (2) of the definition holds as well as the proof that the variety is
non-trivial rests largely on the following claim.

\begin{claim}
For all $m$ there is an element of ${\bf K }$ so that in that algebra
$c_{00} = d_m$ (with the $d_n$'s distinct --- see definition of ${\bf K}$,
$K_0$ )
\end{claim}

The proof of the claim is quite easy. We inductively define an
equivalence relation $\equiv$ on the constants and an interpretation of
the functions as projections. At stage $n$, we define an equivalence
relation $\equiv_n$ on $\{ c_{km},d_m\colon m,k <\omega \}$, so that
$\equiv_n $ is a subset of $\equiv_{n+1}$, all but finitely many equivalence
classes of $\equiv_n$ are singletons and $\Sigma \{ {\rm card}( A ) -1
\colon A \mbox{ an} \equiv_n \mbox{-equivalence class} \} \le 2n$. Moreover
we demand that for each $l<\omega$ no two distinct members of
$\{c_{l,m}:m<\omega\}\cup\{d_m:m<l\}$ are $\equiv_n$-equivalent. To begin we
set $c_{0,0} \equiv_0 d_m$. At stage $n$, there are two possibilities,
either $t_n$ has already been set equivalent to one of $\{d_k \colon k < m_n
\} \cup \{c_{m_n,i} \colon i < k_n\}$ or not. In the first case our
assumption on the arity guarantees that we can make $f_n$ a projection
function and we put $\equiv_n=\equiv_{n-1}$. In the second case there is
some element in $\{c_{m_n,j}:j<k_n\}$ which has not been set equivalent to
any other element. In this case we choose such an element, set it equivalent
to $t_n$ and to some $c_{0,m}$ (for a suitable $m$) and let $f_n$ be the
appropriate projection. In the end let $\equiv$ be \ \ \ \ $x\equiv y\mbox{
\iff \ } (\exists n)(x\equiv_n y)$.\ \ \ \ Let $M$ be the $\tau'$-algebra
with the set of elements $\{c_{n,m},d_n:n,m<\omega\}/\equiv$, functions
$f_n$ as chosen above (note that $f_n$ respects $\equiv$ as it is a
projection) and $c_{n,m},d_n$ interpreted naturally. Note that by the
equation $(*)$ for every $m$, $\{c_{m,n}:n<\omega\}\cup\{d_n:n<m\}$ lists
the members.

It remains to verify that condition (2) is satisfied. Let $I$ denote the
free algebra generated by $\{c_{0,n} \colon n < \go \}$.  Suppose that (2)
is not true. Choose elements $\{e_n \colon n <\go\}$, so that $I \ast
F(\aleph_0)$ is freely generated by the $d_n$'s and the $e_n$'s. So there is
some $m$, so that $c_{0,0}$ is in the subalgebra generated by $\{ d_n \colon
n < m\} \cup \{e_n \colon n <\go\}$. But if we turn to the projection
algebra where $c_{0,0} = d_m$, we have a contradiction (see definition of
${\bf K}$).

Finally we need to see that our variety does not satisfy \cp{2}. We will
prove the following claim which not only establishes the desired result
but shows the limit of our method.

\begin{claim}
If ${\cal V}$ is a variety which is generated by projection algebras then
${\cal V}$ does not satisfy \cp{2}.
\end{claim}

Suppose to the contrary that we had such a variety. Let $I$ and
\(\{x_{in} \colon i <2, n < \go \}\) be an example of \cp{2}.  Choose
\(\{y_n \colon n < \go \}\) so that 
\[\{x_{1n} \colon n < \go\} \cup \{y_n \colon n < \go \}\] 
is a set of free generators for $I$ (we rename $\{h_n:n\in s_i\}$ as
$\{x_{in}:n<\omega\}$.  Notice that if $\gth$ is any verbal congruence (see
below) on $I$ which does not identify all elements then the image of a set
of free generators of a subalgebra will freely generate their image in the
subvariety ${\cal V} / \gth $ defined by the law. (A {\em verbal congruence}
is a congruence which is defined by adding new laws to the variety.) Fix a
vocabulary.

We will show by induction on the complexity of terms $\tau$ that 
\begin{description}
\item[$\otimes$\ \ \ ] if $\cal V$ is a variety generated by projection
algebras and if $X \cup Y$ are free generators of an algebra $A \in {\cal V
}$ , $ a = \tau (\ldots,x_i,\ldots,y_j,\ldots)$, $x_i \in X$, $y_j \in Y$,
$a \in A$ and $X\cup \{a\}$ freely generates a subalgebra of $A$ then $a$ is
in the subalgebra generated by $Y$.
\end{description}

The base case of the induction is trivial. Suppose that $a = f(t_0,\ldots,
t_n)$.  For $i \leq n$ let $\gth_i$ be the congruence on $A$ generated by
adding the law $f(z_0,\ldots, z_n) = z_i$ and let ${\cal V}_i$ be the
subvariety satisfying this law. Since $\cal V$ is generated by projection
algebras so is ${\cal V}_i$, for all $i$.  Furthermore $\cal V$ is the join
of these varieties. In $A/\gth_i$, $a/\gth_i = t_i/\gth_i$. By the inductive
hypothesis we can choose for each $i$, a term $s_i$ using only the variables
from $Y$ so that $A/\gth_i$ satisfies that $t_i/\gth_i = s_i/\gth_i$. Hence
each variety, ${\cal V}_i$ satisfies the law $f(t_0, \ldots, t_n) = f(s_0,
\ldots, s_n)$.  So ${\cal V } $ satisfies the law as well. We have shown
that $a = f(s_0,\ldots, s_n)$, i.e., that $a$ is in the subalgebra generated
by $Y$. So $\otimes$ holds.

Applying the last claim we have, \(\{x_{0n} \colon n < \go\}\) is
contained in the subalgebra generated by \(\{y_n \colon n < \go \}\).
Call the latter subalgebra $B$. Let $F$ denote a countably generated free
algebra. Since $B \ast F$ is isomorphic over $B$ and hence over $\{x_{0n}
\colon n < \go \}$ to $I$, and $\{x_{0n}:n<\omega\}$ generates a free factor
of $I$, necessarily $\{x_{0n} \colon n < \go \}$ freely generates a free
factor of $B \ast F$. Hence \(\{x_{in} \colon i <2, n < \go \}\) freely
generates a free factor of $I \ast F$. Thus we have arrived at a
contradiction. $\Box$
\medskip

\section{Miscellaneous}
\setcounter{theorem}{0}
\setcounter{defi}{0}
One natural question is for which cardinals $\gk$ is every $\gk$-free
algebra of cardinality $\gk$ free (no matter what the variety). By the
singular compactness theorem (\cite{Sing}) every singular cardinal is such a
cardinal. As well it is known that if $\gk$ is a weakly compact cardinal
then every $\gk$-free algebra of cardinality $\gk$ is free.  Some proofs of
this fact use the fact that for weakly compact cardinals we can have many
stationary sets reflecting in the same regular cardinal (see the proof of
\cite{EMB} IV.3.2 for example). It turns out that we only need to have
single stationary sets reflecting. We say that a stationary subset $E$ of a
cardinal $\gk$ {\em reflects} if there is some limit ordinal $\ga < \gk$ so
that $E \cap \ga$ is stationary in $\ga$.  The relevance of the following
theorem comes from the fact that the consistency strength of a regular
cardinal such that every stationary set reflects in a regular cardinal is
strictly less than that of a weakly compact cardinal \cite{refl}. So the
consistency strength of a regular cardinal $\gk$ so that every almost free
algebra of cardinality $\gk$ is free is strictly less than that of the
existence of a weakly compact cardinal.

We separate out the following lemma which will be useful in more than one
setting.

\begin{lemma}
\label{trans-claim}
Suppose $F$ is a free algebra and $G \se H$ are such that $H$ is a free
factor of $F$ and there are $A$, $B$ free subalgebras of $H$ so that $G
\se A$, ${\rm rank}(B) = {\rm card} (H)+\aleph_0$ and $A * B = H$. Then
$G$ is a free factor of $F$ \iff\ $G$ is a free factor of $H$.
\end{lemma}

\proof Obviously if $G$ is a free factor of $H$ then $G$ is a free
factor of $F$. Suppose now that $G$ is a free factor of $F$.  Since $G$ is a
free factor of $F$ we can choose $B_1$ so that $|B| = |B_1|$ and $G$ is a
free factor of $A\ast B_1$. But since $H$ is isomorphic over $A$ to $A \ast
B_1$, $G$ is also a free factor of $H$. \fin

Notice that in the hypothesis of the last lemma the existence of $A$ and $B$
is guaranteed if we assume that $|G| < |H|$.
   
In some varieties a union of an increasing chain of cofinality at least
$\gk$ of $\gk$-pure subalgebras is $\gk$-pure as well. In general
varieties this statement may not be true, in our later work we will need
the following weaker result.

\begin{theorem}
\label{refl-thm}
Suppose $\gk$ is an inaccessible cardinal and $E$ is a subset of $\gk$ such
that every stationary subset of $E$ reflects in a regular cardinal.  If $A$
is an almost free algebra of cardinality $\gk$ and there is a
$\gk$-filtration $(A_\ga\colon \ga < \gk)$ of $A$ such that for all $\ga
\notin E$, $A_\ga$ is $\gk$-pure then $A$ is free.
\end{theorem}

\proof Since $\gk$ is a limit cardinal, $A$ is also $\mbox{L}_{\infty
\gk}$-free. Hence wlog $E$ is a set of limit ordinals and each $A_\ga$ is
free (here we use that $A$ is ``almost free"); for all $\ga < \gb$, $|A_\ga|
< |A_\gb|$ and for all $\ga < \gb$, $A_{\ga+1}$ is a free factor of
$A_{\gb+1}$ and $\alpha\in\lambda\setminus E$ implies that $A_\alpha$ is
$\kappa$-pure. Also wlog 
\begin{description}
\item[$(*)_1$\ \ \ ] if $A_\alpha$ is $\kappa$-pure in $A$ then
$A_{\alpha+1}/A_\alpha$ is free
\item[$(*)_2$\ \ \ ] if $A_\alpha$ is not $\kappa$-pure in $A$ then
$A_{\alpha+1}/A_\alpha$ is not essentially free.
\end{description}
Assume that $A$ is not free. We will use the fact that
\begin{quote}
If $\{A_\ga\colon \ga < \gk\}$ is a filtration of $A$ and $B$ is a
$\gk$-pure subalgebra of cardinality less than $\gk$ then for a club $C$
of $\gk$, $\gb \in C$, cf$\gb=\go$ implies that $A_\gb/ B$ is free.
\end{quote}
Let 
$$E^* = \set{\ga<\kappa}{A_\ga \mbox{ is not a free factor of }A_{\ga
+1}},$$ 
\[C^*=\{\alpha<\kappa:\alpha \mbox{ is a limit cardinal and } \beta<\alpha
\Rightarrow |A_\alpha|<\alpha\}.\]
Now $C^*$ is a club of $\kappa$ and $E^*$ is a stationary subset of $\kappa$
(otherwise $A$ is free). Choose $\gl$ a regular cardinal so that $|A_\gl| =
\gl$ and $(E^*\cap C^*)\cap \gl$ is stationary in $\gl$. If $A_\gl$ is free
then we can find a strictly increasing continuous sequence $\langle
\ga_i\colon i < \gl\rangle$ such that $i < j$ implies $A_{\ga_j}/A_{\ga_i}$
is free. Let $C=\{i <\gl\colon \ga_i =i\}$. Since $(E^*\cap C^*)\cap \gl$ is
stationary there is $\gb\in E^*\cap C^*\cap \gl \cap C$. So we can find
$\gb_n$ for $n < \go$ such that $\gb_n \in C$, $\gb_0 = \gb$ and $\gb_n
<\gb_{n+1}$.

Let $\gb_\go = \cup \{\gb_n \colon n < \go\}$. Then $A_{\gb_\go}/ A_\gb$ is
free by the choice of $C$ (and of the $\alpha_i$'s). Also
$A_{\gb_{n+1}}/A_{\gb_0 +1}$ is free.  Together $A_{\gb_0+1}/ A_{\gb_0}$ is
essentially free, so by $(*)_2$ we know $A_{\gb_0 +1}/ A_{\gb_0}$ is free
which contradicts our choice of $E^*$ and $\beta_0\in E^*$. Hence $A_\gl$ is
not free. \fin

\section{Getting \cp{n}}
\setcounter{theorem}{0}
\setcounter{defi}{0}

In this section we will deal with the problem of deducing \cp{n} from the
existence of a $\gk$-free algebra. We will need to deal with subalgebras of
free algebras.

To handle certain technical details in this section we will deal with
varieties in uncountable languages.  Most things we have done so far
transfer to this new situation if we replace of cardinality $\gk$ by
$\gk$-generated. One trick we will use is to pass from a pair of algebras
$B/A$ to a new algebra $B^*$ by making the elements of $A$ constants in the
new variety. Recall that the notation $B/A$ implies that $B$ and $A$ are
subalgebras of an algebra $C$ and both $A$ and the subalgebra generated by
$B\cup A$ are free.

\begin{defi} Suppose $\cal V$ is a variety and $A$ is a free algebra in
the sense of $\cal V$. Let ${\cal V}_A$ denote the variety where we add
constants for the elements of $A$ and the equational diagram of $A$.
\end{defi}

Notice that any element of ${\cal V}_A$ contains a homomorphic image of $A$.
If $B$ is any algebra which contains $A$ we can view $B$ as a ${\cal V}_A$
algebra.

\begin{proposition}
Suppose $A$ is a $\cal V$-free algebra and $B$ is an algebra which contains
$A$. Then for all $\gk$, $B$ is $\gk$-free in ${\cal V}_A$
\iff\ $B/A$ is $\gk$-free in ${\cal V}$.
\end{proposition}
  
The following lemma is easy and lists some facts we will need. 

\begin{lemma}
Suppose $A \subseteq B$ and both $A$ and $B$ are free algebras on $\gk$
generators. Then the following are equivalent.

(i) every subset of $A$ of cardinality $< \gk$ is contained in a subalgebra
$C$ which is a free factor of both $A$ and $B$.

(ii) every free factor of $A$ which is $\ltc\gk$-generated is also a free
factor of $B$.
\end{lemma}

\proof That (ii) implies (i) is obvious. Assume now that (i) holds
and that $C$ is a free factor of $A$ which is $\ltc\gk$-generated. Let
$D$ be a $\ltc\gk$-generated free factor of both $A$ and $B$ which
contains $C$. Since $A \cong_D B$, $A \cong_C B$. So $C$ is a free factor
of $B$. \fin

If $A$ and $B$ are free algebras which satisfy either (i) or (ii) above, we
will write $A \isub{\gk} B$. This notation is justified since for free
algebras these conditions are equivalent to saying $A$ is an
$\mbox{L}_{\infty \gk}$-subalgebra of $B$. It is possible to give a simpler
characterization of \cp{n}.

\begin{theorem}
\label{cpcrit}
For any variety of algebras, \cp{n} is equivalent to the following
statement. There are countable rank free algebras, $A \isub{\go} B$ and
countable rank free algebras $A_k$ $(k < n)$ so that 
\begin{description}
\item[(i)\ \ \ ] $A = \ast_{k < n} A_k$ and for all $m$, $\ast_{k \neq m}
A_k$ is a free factor of $B$ 
\item[(ii)\ \ ] if $F$ is a countable rank free algebra, then $A$ is not a
free factor of $B\ast F$ (alternatively, $B/A$ is \enf).
\end{description}
\end{theorem}

\proof \cp{n} clearly implies the statement above. Assume that $A$,
$B$, $A_k$ $(k < n)$ are as above. We will show that $B \ast F$ together
with $A$ satisfies \cp{n} with $A,B,B*F$ here corresponding to to $H,I,L$
there.  It is enough to prove that (i) in the statement of \cp{n} holds.  It
suffices to show for all $m$ that if $C$ is a finite rank free factor of
$A_m$ then $\ast_{k \neq m} A_k \ast C$ is a free factor of $B \ast F$.
Choose $Y$ a complementary factor in $B$ for $\ast_{k \neq m} A_k$.  Choose
now finite rank free factors $D$ and $E$ of $\ast_{k \neq m} A_k$ and $Y$
respectively so that $C$ is contained in $D \ast E$ and is a free factor of
$D \ast C$. Clearly $D\ast E$ is a free factor of $B$ and also of $D\ast Y$,
hence we have $B \ast F \cong_{D \ast E} D \ast Y
\ast F$. Now $C \ast D$ is a free factor of $A$, so as $A$ is an $L_{\infty
, \omega }$-submodel of $B \ast F$, all are countable generated, clearly
$C\ast D$ is a free factor of $B\ast F$. By the last two sentences (as
$C\ast D\subseteq D\ast E$) we have that $C\ast D$ is a free factor of $D
\ast Y \ast F$. Also $\ast_{l\neq m}A_l$, $D\ast Y\ast F$ are freely
amalgamated over $D$, $D$ is a free factor of both and $D\subseteq D\ast
C\subseteq D\ast Y\ast F$ are free. $D$ is a free factor of $D\ast C$,
$D\ast C$ is a free factor of $D\ast Y\ast F$; together $\ast_{l\neq m}\ast
C$ is a free factor of $B\ast F$. So we have finished. \fin

We next have to consider pairs (and tuples). The following two facts are
standard and proved analogously to the results for algebras (rather than
pairs).


\begin{lemma}
Suppose $B/A$ is $\gk^+$-free. Then it is strongly $\gk$-free (i.e.
\lsub{\gk}-equivalent to a free algebra).  
\end{lemma}

\begin{corollary}
\label{isubchain}
Suppose $\gk$ is regular and $A \se B$. If $B/A$ is $\gk$-free and $|B| =
\gk$ then $B = \cup_{\ga < \gk} B_{\ga}$ (continuous) where $A = B_0$,
$B_{\ga + 1}$ is countably generated over $B_\ga$ and for all $\ga$,
$B_\ga \isub{\go} B$.
\end{corollary}

We now want to go from the existence of certain pairs to \cp{n} for various
$n$. The difficulty is in suitably framing the induction hypothesis. We
define the pair $B/A$ to be $\ha_0$-free if $A \isub{\go} A + B$. In order
to state our result exactly we will make an {\em ad hoc} definition.

\begin{defi} We say $\gk$ {\em implies \cp{ n,m}} if: $\gk$ is regular,
and for any variety ${\cal V}$ if $(*)_{\gk , m}$ below holds then the
variety satisfies \cp{n+m} where:
\begin{description}
\item[$(*)_{\gk,m}$\ \ \ ] there are free algebras (free here means 
in ${\cal V }$ ) $A,B , F_0, F_1,F_2,\ldots ,F_m$ such that
\begin{description}
\item[(a)\ ] all are  free 
\item[(b)\ ] all have dimension $\gk$ ( i.e., a basis of cardinality $\gk$)
\item[(c)\ ] $A$ is a subalgebra of  $B$
\item[(d)\ ] $A$ is the free product of   $F_0,F_1,F_2,\ldots, F_m$
\item[(e)\ ] if $m>0$, for  $k \in \{ 1,2,\ldots,m \}$, $ B$ is free
over the free product of $\{ F_i : i \le m \mbox{ but } i \mbox{ is not
equal to }k \}$
\item[(f)\ ] $B/A$  is $\gk$-free but not essentially free.
\end{description}
\end{description}

We say $\gk$  {\em implies \cp{n}} if for every $m$  it implies
\cp{n,m}. 
\end{defi}

Remark: {\em 1)} We can weaken the demand on the $F_i$ to having the
dimension be infinite. Note, as well, that there is no demand that the free
product of $\{F_i\colon i \neq 0\}$ is a free factor of $B$.\\
{\em 2)} Remember that if ${\cal V}$ satisfies \cp{n+1} then ${\cal V}$
satisfies \cp{n}.

\begin{proposition} 
$\aleph_0$  implies  \cp{0}.
\label{base}
\end{proposition}

\proof Without loss of generality we can assume that $B$ is
isomorphic to $B*F$ over $A$ where $F$ is a countable rank free algebra.
There are two cases to consider. First assume that the free product of
$\{F_i\colon i \neq 0\}$ is a free factor of $B$. In which case by
Theorem~\ref{cpcrit} we have an example of \cp{m+1} (hence ${\cal V}$
satisfies \cp{m}). Next assume that the free product of $\{F_i\colon i \neq
0\}$ is a not free factor of $B$. We claim that $A^*$ and $B$ are an example
of \cp{m}, where $A^*$ is the free product of $\{F_i\colon i \geq 1\}$.  All
that we have to check (by \ref{cpcrit}) is that for all $k$, such that $ 1
\leq k\leq m$, the free product of $\{F_i\colon i \geq 1, i \neq k\}$ is a
free factor of $B$. But this is part of the hypothesis. \fin 

Note that in the definition we can allow to increase all dimensions to be
just at least $ \gk$ except that $B$ should be generated by $A$ together
with a set of cardinality $\gk$.

We will take elementary submodels of various set-theoretic universes and
intersect them with an algebra.

\begin{proposition}  
\label{enf-prop}
Suppose that $A$ and $B$ are free algebras and $B/A$ is \enf\ and $B$ is
$\gk$-generated over $A$. If $N \prec (\He{\chi})$, where $A, B \in N$, $\gk
+ 1 \se N$ and $\chi \geq |A| + \gk$, then $(B \cap N)/ (A \cap N)$ is \enf.
Furthermore if $B/A$ is $\gk$-free, then so is $(B \cap N)/ (A
\cap N)$.
\end{proposition}

\proof First deal with the first assertion. Let $Y \in N$ be a free basis of
$A$. (Note that such a $Y$ must exist since $A \in N$.) So $A \cap N$ is
freely generated by $Y\cap N$. Without loss of generality, we can assume
that $B$ is isomorphic over $A$ to $B \ast F$ where $F$ is a free algebra of
rank $\gk$. Under this assumption, $(B \cap N)/(A \cap N)$ is \enf\ \iff\ it
is not free.  Suppose that $(B \cap N)/(A \cap N)$ is free. Then we can find
$Z\subseteq(N\cap B)$ so that $Z \cup (Y \cap N)$ is a basis of $N \cap B$.

We first claim that $Z \cup Y$ is a free basis for the algebra it generates.
If not, then there is some finite $Z_1 \se Z$ and finite $Y_1\se Y$ so that
$Z_1 \cup Y_1$ is not a free basis for the algebra it generates. That is it
satisfies some equation which is not a law of the variety. By
elementariness, we can find $Y_2\se N \cap Y$ so that $Z_1\cup Y_2$
satisfies the same law. This is a contradiction. 

To finish the proof we must see that $Z \cup Y$ generates $B$.  Choose $X
\in N$ of cardinality $\gk$ so that $X \cup A$ generates $B$. Since $ \gk
+ 1 \se N$, $X \se N$. Hence $X \se (B \cap N)$. As $Z\cup(Y\cap N)$
generates $B\cap N$ and $Y$ generates $A$ we conclude that $X\cup A$ (and
hence $B$) is contained in the algebra generated by $Z\cup Y$. 

The second statement is very simple to prove. If $B/A$ is $\gk$-free, choose
in $N$, $X$ of cardinality $\gk$ so that $B$ is generated over $A$ by $X$
and a sequence $(X_\ga \colon \ga < \gk)$ which witnesses that $B/A$ is
$\gk$-free. Since each $X_\ga \se N$, ${(B \cap N)/ (A \cap N)}$ is
$\gk$-free. \fin

\begin{theorem}
\label{add1-thm}
(1) $\aleph_1$ implies \cp{1}.

\noindent
(2) Suppose $\gk$ is a regular cardinal. Assume that for every variety
${\cal V}$ and free algebras $A,B$ in ${\cal V}$, if $B/A$ is $\gk$-free
\enf\ of cardinality $\gk$ then there are: $\gl$ which implies
\cp{n,m+1}, $\chi$ large enough and $M \prec (\He{\chi}, <)$ of
cardinality $<\kappa$, $M\cap\kappa$ an ordinal $\geq\lambda$ so that $A, B
\in M$ and $B\cap M/A$ is $\lambda$-pure in $B/A$ and there is an elementary
submodel $N\prec(\He{\chi},<)$ such that $A,B,M\in N$, $\lambda+1\subseteq
N$, $|N|=\lambda$ and $N\cap B/A\cup(B\cap N)$ is almost free \enf\ .  Then
$\gk$ implies \cp{n+1,m}.

\noindent
(3)  Suppose $\gk$ is a regular cardinal. Assume that for every variety
${\cal V}$ and free algebras $A,B$ in ${\cal V}$, if $B/A$ is $\gk$-free
\enf\ of cardinality $\gk$ then there are: $\gl$ which implies \cp{n}, $\chi$
large enough and $M \prec (\He{\chi}, <)$ of cardinality $<\gk$,
$M\cap\kappa$ an ordinal $\geq\lambda$, so that $A, B \in M$ and $B\cap M/A$
is $\gl$-pure in $B/A$ and there is an elementary submodel
$N\prec(\He{\chi},<)$ such that $A,B,M\in N$, $\lambda+1\subseteq N$,
$|N|=\lambda$ and $N\cap B/A\cup(B\cap N)$ is almost free \enf\ .  Then
$\gk$ implies \cp{n+1}.
\end{theorem}

\proof By Corollary~\ref{isubchain} and Proposition~\ref{base}, (1)
is a special case of (2).  Also part (3) follows from part (2), by the
definitions. So we will concentrate on that case.

Consider an instance of checking that $\gk$ implies \cp{n+1,m}, i.e., we are
given $A,B,F_0,F_1,\ldots,F_m$ as in the definition of $(*)_{\kappa,m}$.
Let $M,\chi ,\gl$ be as guaranteed in the assumption of the theorem and let
$N$ be an elementary submodel of the $\chi$ approximation to set theory to
which $A,B, M,F_0,F_1,\ldots,F_m$ belong, $N$ has cardinality $\gl$, 
$\lambda+1\subseteq N$ and $(N \cap B ) / ( A \cup (B \cap M ) )$ is almost
free essentially non free of dimension $\gl$.

Note first that there is a filtration $(B_i \colon i < \gk)$ of $B$ such
that for all $i$, $B_i/A$ is free. We assume that the filtration is in
$M$, so there is some $i$ such that $M \cap B = B_i$. So in particular, $M
\cap B/A$ is free. It is now easy to see that the algebra generated by
$(M\cap B) \cup A$ is the free product over $A \cap M$ of $A$ and $M\cap B$.
More exactly it suffices to show that there any relation between elements is
implied by the laws of the variety and what happens in $A$ and $B\cap M$.
Fix $Y\in M$ a set of free generators of $A$.  As we have pointed out before
$Y \cap M$ freely generates $A \cap M$.  As well, since $B\cap M/A$ is free,
for any finite set $C \se B \cap M$, there is a finite subset $D \se B\cap
M$ so that $C\subseteq D+(A\cap M)$ and $D/A\cap M$ is free. Let $Z\se M$ be
such that $Z\cup (Y \cap M)$ is a set of free generators for $A + D$. To
finish the proof that $A + B\cap M$ is the free product of $A$ and $B\cap M$
over $A\cap M$, it suffices to see that $Y \cup Z$ freely generates $A + D$.
The set obviously generates $A + D$. By way of obtaining a contradiction
assume that a forbidden relation holds among some elements of $Z$ and some
elements of $Y$. Then since $M$ is an elementary submodel of an
approximation to set theory, the elements of $Y$ can be taken to be in $Y
\cap M$, which contradicts the choice of $Z$.  Finally note that since
$B\cap M / A$ is free, $B\cap M/A\cap M$ is essentially free.

Let $A_0$ be $A \cap M$ and let $B_0$ be $B \cap M$.  Let $A_1$ be the
subalgebra of $B$ generated by $A \cup B_0$.

As each $F_k$ ($k \le m$) is free and we can assume belongs to $M$, clearly
$$F^0_k =^{df} F_k \cap M$$ is free of dimension $\gl$ and $F_k$ is the free
product of $F^0_k$ and some free $F^1_k$ which has dimension $\gk$.

Without loss of generality, each $F^1_k$ belongs to $N$. Let $B_1$ be the
subalgebra of $B$ generated by $B_0 \cup F^1_0 $. Since $B_0 / A_0$ is
essentially free, $B_1/A_0$ is free. Let $F^1_{m+1}$ be a free subalgebra of
$B_1$ such that $B_1$ is the free product of $A_0$ and $F^1_{m+1}$.

Without loss of generality, $F^1_{m+1} \in N$.  For $k \le m+1$, let $F^*_k$
be $F^1_k \cap N$ if $k > 0$ and $A_0$ if $k = 0$. Let $B^*$ be $B \cap N$
and let $A^*$ be the subalgebra of $B$ generated by $[A \cup ( B \cap M
)]\cap N$.

That $B^*,A^*,F^*_k $ (for $k \le m+1$) are free should be clear, as well as
the fact that $A^*$ is the free product of $\{F^*_k\colon k
\leq m+1\}$. Furthermore $B^*/A^*$ is almost free but \enf. It
remains to prove that if $k \le m+1$ is not zero, then $B^*$ is free
over the free product of $\{F^*_i\colon i \neq k\}$. After we have
established this fact we can use that $\gl$ implies \cp{n, m+1}, to
deduce that \cp{n+m+1}\ holds in the variety. 

Assume first that $k \le m $. We know that $B$ is free over the free product
of $\{F_i\colon i \neq k\}$. So $B \cap N$ is free over the algebra
generated by $(B \cap M) \cup \cup \{F_i \cap N : i \le m, i \neq k \}$. But
the algebra generated by $(B\cap M) \cup *_{i \leq m, i\neq k}F_i$ is the
same as $A_0 * F_{m+1}^1**_{0 <i \leq m, i\neq k}F_i^1$.  Next assume that
$k=m+1$. We must show that $B \cap N$ is free over $A_0 * F^*_1 * \cdots *
F^*_m$. As above $B\cap N$ is free over $A$ and so as before $B\cap N$ is
essentially free over $A \cap N$. Since $A \cap N = A_0 * F^*_1 * \cdots *
F^*_m * (F^1_0 \cap N)$, $B\cap N$ is essentially free over $A_0 * F^*_1 *
\cdots * F^*_m$. So replacing, if necessary, $B\cap N$ by its product with a
free algebra we are done.
\fin

\noindent{\sc Discussion:}\ \ \ In order to get a universe where the
existence of an $\aleph_n$-free\enf\ algebra implies \cp{n}, we will use
various reflection principles.  We will consider sentences of the form $Q_1
X_1 Q_2 X_2 \ldots Q_n X_n\psi(X_1,\ldots X_n)$, where $Q_1, \ldots Q_n$ are
either $aa$ or $stat$ and $\psi(X_1, \ldots, X_n)$ is any (second-order)
sentence about $ X_1,\ldots, X_n$ (i.e., $\psi$ is just a first order
sentence where $X_1,\ldots, X_n$ are extra predicates). We call this
language $L_2(aa)$.  To specify the semantics of this language we first fix
a cardinal $\gl$, and say in the $\gl$-interpretation, a model $A$ satisfies
$aa\, X\, \psi(X)$ if there is a closed unbounded set, $\cal C$ in ${\cal
P}_{<\gl}(A)$ so that $\psi(X)$ for all $X \in \cal C$. Similarly $stat$
means ``for a stationary set''. To denote the $\gl$ interpretation we will
write $L_2(aa^\gl)$ Notice that the L\'evy collapse to $\gl^+$ preserves any
statement in the $\gl$-interpretation.

\begin{defi}
 For regular cardinals $\gk , \gl$, let $\gD_{\gk\gl}$ denote the following
principle:
\begin{quote}
Suppose $A$ is a structure of whose underlying set is $\gk$ and $\gf$ is any
$L_2(aa^\gl)$-sentence. Let $C$ be any club subset of $\gk$. If $A$
satisfies $\gf$ then there is a substructure in $C$ of cardinality $\lambda$
which satisfies $\gf$. 

\end{quote}
\end{defi}

This principle is adapted from the one with the same name in
\cite{MaSh}. They use it to show that that $\ha_{\go^2 +1}$ may be
outside the incompactness spectrum of abelian groups. 

\begin{theorem}[\cite{MaSh}] 
(1) If $\gk$ is weakly compact and $\gk$ is collapsed (by the L\'evy
collapse) to $\aleph_2$, then $\gD_{\aleph_2\aleph_1}$ holds.  

(2) Suppose it is consistent that that there are $\aleph_0$ supercompact
cardinals.  Then it is consistent that for every $m > n$, $\gD_{\aleph_{m},
\aleph_n}$ holds.
\end{theorem}

\begin{theorem}
(1) Suppose $\gD_{\aleph_2\aleph_1}$ holds. If in a variety there is an
$\aleph_2$-free \enf\ algebra of cardinality $\aleph_2$, then \cp{2} holds.

(2) Suppose $m > 1$ and for every $m \geq n > k$, $\gD_{\aleph_{n},
\aleph_k}$  holds. If in a variety there is an $\aleph_m$-free \enf\ algebra
of cardinality $\aleph_m$, then \cp{m} holds. 
\end{theorem}

\proof (1) Suppose $A$ is an $\aleph_2$-free \enf\ algebra of
cardinality $\aleph_2$. Without loss of generality we can assume that $A
\cong A \ast F$ where $F$ is a free algebra of cardinality $\aleph_2$. Hence
if $A^*\subseteq A$, $|A^*|<|A|$, $\mu\leq\aleph_2$ then $A/A^*$ is
essentially $\mu$-free \iff\ $A/A^*$ is $\mu$-free.  Consider $C =
\set{X}{|X| = \aleph_0 \mbox{ and } A/X\mbox{ is } \aleph_1\mbox{-free}}$.
We claim that $C$ must contain a club. Otherwise we would be able to reflect
to a free subalgebra of cardinality $\aleph_1$ which satisfies (in the
$\aleph_0$ interpretation)
\[ stat X stat Y (Y/X \mbox{ is \enf}) \]
which contradicts the freeness of the subalgebra.

Consider now a filtration $\set {A_\ga}{\ga < \go_2}$ (such that each
$A_\alpha$ is free). If
\[\set{\ga}{A/A_\ga \mbox{ is } \aleph_1\mbox{-free and \enf}}\] 
contains a club, hence $\{\alpha: A/A_\alpha\mbox{ is } \aleph_1$-free and
not almost free$\}$ is stationary and then we are done since in this case by
Theorem~\ref{add1-thm}, $\ha_1$ implies \cp{1} and we have have an instance
of $(*)_{\ha_1, 1}$. Hence the variety satisfies \cp{2}.

Assume the set does not contain a club, i.e.,
\[S = \set{\ga}{A/A_\ga \mbox{ is essentially not }\aleph_1\mbox{-free}}\]  
is stationary. Let $\open P$ be the L\'evy collapse of $\aleph_2$ to
$\aleph_1$.  Then, by the previous paragraph,
\[\force_{{\open P}} A\mbox{ is free}.\]
Let $\dot{X}$ be the name for a free basis of $A$. Choose some cardinal
$\chi$ which is large enough for $\open P$. Let $M \prec (\mbox{H}(\chi),
\in, <)$ where $<$ is a well-ordering, everything relevant is an element
of $M$, $\go_1 \se M$, and $M \cap\go_2 = \gd \in S$.  Next choose a
countable $N \prec (\mbox{H}(\chi), \in, <)$ so that $M \in N$. Hence $(N
\cap A)/A_\delta$ is \enf. We will contradict this statement and so
finish the proof. By Proposition~\ref{enf-prop}, $(N \cap A)/(N\cap
A_\gd)$ is \enf. 

Since $N \cap A \in C$, we have $A/(N \cap A)$ is $\aleph_1$-free. We shall
show that $A/(N \cap A_\gd)$ is $\ha_1$-free. Let us see why this finishes
the proof. By the two facts we can find a countable subalgebra $B$ so that
$B/(N \cap A)$ and $B/(N \cap A_\gd)$ are both free. But since $B = (N \cap
A) \ast F$ for some free algebra $F$, we would contradict the fact $(N \cap
A)/ (N \cap A_\gd)$ is \enf.

Let $p$ be an $N \cap M$-generic condition. Then
\[p \force N \cap M \cap A\ (= N \cap A_\gd) \mbox{ is generated by }
N \cap A_\gd \cap \tilde{X}. \] So $p$ forces that $A/(N \cap A_\gd)$ is
$\aleph_1$-free. But being $\aleph_1$-free is absolute for L\'evy
forcing.

The proof of part (2) is similar. \fin

The situation in part (2) of the theorem above is perhaps the most
satisfying. On the other hand we need very strong large cardinal assumptions
to make it true. (It is not only our proof which required the large
cardinals but the result itself, since if the conclusion of (2) is satisfied
then we have for all $m$, any stationary subset of $\aleph_{m+1}$ consisting
of ordinals of cofinality less than $\aleph_m$ reflects.) It is of interest
to know if the classes can be separated via a large cardinal notion which is
consistent with V = L.  Rather than stating a large cardinal hypothesis we
will state the consequence which we will use.

\begin{defi} 
Say a cardinal $\gm$ is an {\em $\go$-limit of weakly compacts} if there are
disjoint subsets $S$, $T_n$ $(n < \go)$ of $\gm$ consisting of inaccessible
cardinals so that
\begin{trivlist}
\item[1.] for every $n$ and $\gk \in T_n$ and $X \se \gk$, there is
$\gl \in S$ so that\\ 
$(V_\gl, X \cap \gl, \in) \prec_{\Sigma^1_1} (V_\gk,X,\in)$
\item[2.] for every $n$ and $\gk \in T_{n+1}$, $T_n \cap \gk$  is
stationary in $\gk$.
\end{trivlist}
\end{defi}

Notice that (1) of the definition above implies that every element of $T_n$
is weakly compact.

\begin{theorem}
Suppose that $\gm$ is a $\go$-limit of weakly compacts and that GCH holds.
Let $S$, $T_n$ $(n < \go)$ be as in the definition. Then there is a forcing
extension of the universe satisfying: for all $n$ and $\gk \in T_n$ there is
a $\gk$-free \enf\ algebra of cardinality $\gk$
\iff\ \cp{n+1} holds.
\end{theorem}

\proof The forcing notion will be a reverse Easton forcing of length
$\gm$. That is we will do an iterated forcing with Easton support to get our
poset $\open P$. The iterated forcing up to stage $\ga$ will be denoted
$\open P_\ga$ and the iterate at $\ga$ will be $\dot{Q}_\ga$.  For $\ga$
outside of $S \cup \cup_{n<\go} T_n$, let $\dot{Q}_\ga$ be the ${\open
P}_\ga$-name for the trivial poset. For $\ga \in S$, let $\dot{Q}_\ga$ be
the ${\open P}_\ga$-name for the poset which adds a Cohen generic subset of
$\ga$. For $\ga \in T_0$, let $\dot{Q}_\ga$ be the ${\open P}_\ga$-name for
the poset which adds a stationary non-reflecting subset of $\ga$ consisting
of ordinals of cofinality $\go$. Finally for $\ga \in T_{n+1}$, let
$\dot{Q}_\ga$ be the ${\open P}_\ga$-name for the poset which adds a
stationary non-reflecting subset of $\ga$ consisting of ordinals in $T_n$.
We will refer to this set as $E_\ga$

\noindent The first fact that we will need is essentially due to Silver and
Kunen (see \cite{Ku}).

\noindent {\bf Fact}. 
{\em 1)}\ \ \ For all $n$ and $\gk \in T_n$, if $\dot{Q}$ is the
${\open P}_\gk$-name for the forcing which adds a Cohen subset of $\gk$,
then
\[\force_{{\open P}_\gk \ast \dot{Q}} \gk \mbox{ is weakly compact}.\]
{\em 2)}\ \ \ If $\dot{R}$ is the ${\open P}_\ga \ast \dot{Q}_\gk$-name for
the forcing which shoots a club through the complement of $E_\gk$, then
\[\force_{{\open P}_\gk} \dot{Q}_\gk \ast \dot{R} \mbox{ is equivalent
to } \dot{Q}.\]

We now want to work in the universe $\mbox{V}^{\open P}$. It is easy to to
see that each $E_\ga$ is a stationary non-reflecting set (since the stages
of the iteration after ${\open P}_{\ga +1}$ add no subsets of $\ga$). We
claim that for all $n$ and $\ga \in T_n$, if $D
\se \ga$ is a stationary set and $D$ is disjoint from $E_\ga$ then $D$
reflects in a regular cardinal. This is easy based on the fact. It is
enough to work in $\mbox{V}^{{\open P}_{\ga +1}}$. Let $R$ be the poset
which shoots a club through the complement of $E_\ga$. After forcing with
$R$, $D$ remains stationary. On the other hand, $\ga$ becomes weakly
compact.  So in the extension $D$ reflects and so it must reflect before
we force with $R$.

It is standard (cf.\ Theorem~\ref{existence}) to show that for all $n$ and
$\gk \in T_n$ if \cp{n+1} holds then there is a \enf\ algebra of cardinality
$\gk$ which is $\gk$-free. In fact we can construct such an algebra to have
$E_\gk$ as its $\gG$-invariant. To complete the proof we will show by
induction on $n$ that if $\gk \in T_n$ then $\gk$ implies
\cp{n}.  For $n = 0$, there is nothing to prove. Suppose the result is
true for $n$ and that $B/A$ is $\gk$-free for some $\gk \in T_{n+1}$. By
Theorem~\ref{refl-thm}, $\gG(B/A) \subseteq \dot{E_\gk}$. By the proof of
Theorem~\ref{refl-thm}, we can write $B$ as $\cup_{\ga < \gk} B_\ga$ a
continuous union of free algebras, so that for all $\ga$, $B_{\ga + 1} + A$
is $\gk$-pure and if $\gl$ is a regular uncountable limit cardinal, then
$B_{\gl +1} + A/B_\gl +A$ is $\gl$-free of cardinality $\gl$ and
\enf\ \iff\ it is not free.  By Theorem~\ref{add1-thm} we are done. \fin

That some large cardinal assumption is needed in the previous theorem is
clear. For example, if there is no Mahlo cardinal in L, then every
uncountable regular cardinal has a stationary subset consisting of ordinals
of cofinality $\go$ which does not reflect. So if there is no Mahlo cardinal
in L, then the \enf\ incompactness spectrum of any variety which satisfies
\cp{1} is the class of regular uncountable cardinals. As we shall see in the
next theorem, the existence of a Mahlo cardinal is equiconsistent with the
existence of a cardinal $\gk$ which is in the \enf\ incompactness spectrum
of a variety \iff\ the variety satisfies \cp{2}. However the situation with
\cp{3} and higher principles seems different. It seems that the existence of
a cardinal which is in the \enf\ spectrum of a variety \iff\ it satisfies
\cp{3} implies the consistency of the existence of weakly compact cardinals.

\begin{theorem}
The existence of a Mahlo cardinal is equiconsistent with the existence of a
cardinal $\gk$ which is in the \enf\ incompactness spectrum of a variety
\iff\ the variety satisfies \cp{2}.
\end{theorem}

\proof We can work in L and suppose that $\gk$ is the first Mahlo
cardinal. Fix $E \se \gk$ a set of inaccessible cardinals which does not
reflect. By Theorem 12 of \cite{refl}, there is a forcing notion which
leaves $E$ stationary so that every stationary set disjoint to $E$ reflects
in a regular cardinal. Now we can complete the proof as above to show that
$\gk$ is as demanded by the theorem. \fin


\shlhetal
\end{document}